\documentclass{amsart}
\usepackage{color}
\usepackage{multirow}
\newcommand{\lra}{\longrightarrow}
\newcommand{\Q}{\mathbb Q}
\newcommand{\F}{\mathbb F}

\newcommand{\PP}{\mathbb P^1}
\newcommand{\p}{\mathfrak p}

\newcommand{\rk}{\mbox{rk}}
\newcommand{\MW}{\mbox{MW}}

\newtheorem{thm}{Theorem}
\newtheorem{prop}[thm]{Proposition}
\newtheorem{lem}[thm]{Lemma}
\theoremstyle{remark}
\newtheorem{rem}[thm]{Remark}
\begin{document}

\title{Non-liftable Calabi-Yau spaces}


\author{S\l awomir Cynk}
\address{Institute of Mathematics, Jagiellonian University, ul. \L
ojasiewicza 6,30-348 Krak\'ow,
Poland\newline\indent
Institute of Mathematics of the
Polish Academy of Sciences, ul. \'Sniadeckich 8, 00-956 Warszawa, Poland}
\email{cynk@im.uj.edu.pl}

\author{Matthias Sch\"utt}
\address{Institut f\"ur Algebraische Geometrie, 
Leibniz Universit\"at  Hannover, Welfengarten 1, 30167 Hannover, Germany} 
\email{schuett@math.uni-hannover.de}


\thanks{Funding from MNiSW under grant no N N201 388834 and 
DFG under grant Schu 2266/2-2 is gratefully acknowledged.}

\date{\today}

\begin{abstract}
We construct many new non-liftable three-dimensional Calabi-Yau spaces
in positive characteristic. 
The technique relies on lifting a nodal model to a smooth rigid
Calabi-Yau space over some number field as introduced in \cite{CvS}. 
\end{abstract}

\maketitle

\section{Introduction} 

In recent years,
Calabi-Yau varieties have been studied extensively in the arithmetic context.
Following the modularity of elliptic curves, there has been great
progress on (singular) K3 surfaces  
\cite{ES} and (rigid) Calabi-Yau threefolds.
The latter have two-dimensional $H^3$ which thus is related to
automorphic forms for $\mbox{GL}(2)$ in the number field case.  
If the variety is defined over $\Q$, this reduces to classical
modularity and nowadays follows from Serre's conjecture \cite{KW}. 

It was a recent insight that modularity could be related to
non-liftability at some bad primes. 
After the first results of Hirokado \cite{Hirokado} and Schr\"oer
\cite{Schroeer}, 
there were constructions of non-liftable Calabi-Yau spaces along these lines by Schoen \cite{Schoen2}
and independently by the first author and van Straten \cite{CvS}.
These constructions are based on specific models as double octics or fiber products of specific rational elliptic surfaces.

The aim of this paper is to extend the techniques and constructions
from \cite{CvS} to other fiber products. 
We shall produce a great number of new examples of non-liftable Calabi-Yau spaces.
In brief, our results are summarized as follows:

\begin{thm}
\label{thm}
For each prime $p<100$ with the exception of $53, 83$ there is a
non-liftable three-dimensional Calabi-Yau space in characteristic
$p$. 
The same  holds for all nodal primes in Table~\ref{tab1} and~\ref{tab3}.
\end{thm}

The main idea of the construction is to exhibit a rigid Calabi-Yau
space over some number field that reduces to a nodal model of the
given Calabi-Yau space $X$ in characteristic $p$.  Then the proof of
non-liftability of $X$ follows the lines of \cite{CvS} which we review
in Section~\ref{sec:prel}.  We shall work with fiber products of
elliptic surfaces.  This approach, introduced by Schoen
in~\cite{Schoen}, will be reviewed in Section~\ref{s:setup}.
Sections~\ref{s:4} and~\ref{s:5} introduce the rational elliptic
surfaces that we shall be concerned with.  The proof of
Theorem~\ref{thm} will be completed in section~\ref{s:III}.

As indicated, the rigid Calabi-Yau spaces that we shall construct might also be
interesting from the view-point of modularity.  Namely, many of them
can be defined over $\Q$ and are thus modular.  It came a bit as a
surprise for us how big the bad primes could actually get
(see Tables \ref{tab1}, \ref{tab3}).  Since
these primes have to appear in the level of the corresponding modular
form, this brings back the question whether up to twisting, there are
infinitely many Hecke eigenforms of weight $4$ (or any integer weight greater
than two) with rational coefficients.

\section{Small resolutions in the algebraic setup}
\label{sec:prel}

In this paper we shall construct non--liftable Calabi--Yau spaces as
small resolutions of fiber products of rational elliptic surfaces with
section.  The method was introduced for semistable elliptic surfaces
(all singular fibers of Kodaira type $I_n$) in \cite{Schoen} with the
goal of producing Calabi-Yau threefolds with certain Euler numbers.  In
\cite{Schuett}, twisted fiber products were used to produce new
modular varieties over $\Q$.  In the present paper we want to apply these
techniques, but allow certain types of unstable fibers.  We are
specially interested in rational elliptic surfaces with four or five
singular fibers, one of them of Kodaira type $II$ and the others
semistable.

Let $Y, Y'$ denote elliptic surfaces over $\PP$ with section.  The
fiber product
\[ 
X=Y\times_{\PP} Y'
\] 
need not be smooth; it is singular exactly at points corresponding
to singular points of fibers in both factors.  
In the case of semistable fibers, the fiber
product has type $A_1$ singular points (nodes) as the only
singularities.  A three--dimensional node admits a so--called small
resolution; here the node is replaced by a smooth variety of
codimension two.  In the complex analytic setup a node is locally
analytically isomorphic to the singularity of the quadratic cone, so
it can be resolved by blowing--up an analytic divisor.  A small
resolution of a projective variety in general is not projective.  A
criterion for projectivity of a small resolution of a fiber product
of semistable elliptic surfaces is given in \cite{Schoen}.

In \cite{Artin} Artin gave an algebraic version of the above
construction.  In the category of algebraic spaces (\cite{Artin2}), it
follows from the following theorem:

\begin{thm}[\mbox{\cite[Thm.~4.9]{Artin2}}]
Let $Y'$ be a closed subspace of
    $X'$, and let $\mathfrak f:\mathfrak X'\lra\mathfrak X$ be a
    formal modification, where $\mathfrak X'$ is the formal completion
    of $X'$ along $Y'$. There exists a modification $f:X'\lra X$
    inducing $\mathfrak f$. It is unique up to isomorphism.
\end{thm}

If the fiber product contains a product of two fibers of type $II$, then
it has a singularity of type $D_4$, i.e.\ with a local analytic
equation
\begin{eqnarray}
\label{eq:D4}
x^2+y^2+z^3+t^3=0. 
\end{eqnarray}
As in the nodal case, this singularity admits a small resolution.
This is achieved by first blowing up the local (analytic) divisor
$x+iy=z+t=0$ and then performing a small resolution of the resulting
node. Contrary to the semistable case, this resolution is never
projective.  In the complex analytic category we obtain a compact (but
non-k\"ahler) Moishezon manifold, in the algebraic category we obtain
an algebraic space.

In \cite{Schoen}, Schoen showed that if $Y$ and $Y'$ are rational and
semi-stable, then $X$ is simply connected with trivial dualizing
sheaf.  Hence a small resolution gives a Calabi-Yau threefold (or a
Calabi-Yau space).  This fact was used in \cite{CvS} to produce
non-liftable fiber products thanks to the following key criterion.

\begin{thm}[{\cite[Thm.~4.3, Cor.~4.4]{CvS}}]\label{thm:cvs}
  Let $X$ be a smooth rigid Calabi--Yau space defined over a finite
  ring extension $R$ of $\mathbb Z$. If $\mathfrak p$ is a prime ideal in
  $R$ such that the reduction $X_{\mathfrak
    p}:=X\otimes_R(R/\mathfrak p)$ of $X$ modulo $\p$ has an node as
  the only singularity, then a small resolution of $X_{\mathfrak p}$
  is a Calabi--Yau space that cannot be lifted to characteristic~0.
\end{thm}

\subsection{Behavior under reduction}
\label{ss:359}

Given the small resolution of a fiber product over some number field,
there are two ways for it to attain bad reduction modulo some prime
$\p$.  The first possibility is that one of the factors $Y, Y'$
degenerates mod $\p$ so that the elliptic structure is lost.  This is
usually a fairly rare situation.  
The more common way to attain bad reduction is a
fiber degeneration.  Here fibers in one or both factors are merged
upon reduction.  Often this results in an additional node in positive
characteristic.  However, not every change of structure of the fiber
product necessarily implies an additional node.  In this paragraph, we
comment on this and similar situations.

Consider the situation where one of the surfaces
has fibers of types $I_1$ and $I_2$ at two points
of $\mathbb P^1$ that are merged together to a
type $III$ fiber under the reduction modulo
$\mathfrak p$, whereas the second one  has fibers of types $I_0$ and $I_1$ at these two points. Then on the fiber
product $X$ in characteristic zero we have
to distinguish two cases.

Assume that $X$
has a  fiber of type $I_2\times
I_1$ with two singular points. 
This situation is sketched as follows:
\begin{eqnarray}
\label{eq:former}
\left.\begin{matrix}
I_2\times I_1\\
\\
I_1\times I_0
\end{matrix}
\;\;
\right\}
\;\;
III\times I_1.
\end{eqnarray}
Either singularity admits a (non--projective)
small resolution by blowing--up the product of one
local branch of the $I_1$ singularity and one of the
components of the $I_2$ resp.~$III$ fiber. 
Here the two local analytic
divisors at the two singular points in
characteristic zero are merged together to a
smooth divisor in the fiber product modulo
$\mathfrak p$. 
Blowing--up this smooth divisor
yields the required small resolution. 
Consequently the reduction
of the small resolution of the fiber product is a
small resolution of the reduction, so the reduction
modulo $\mathfrak p$ is smooth.

This should be contrasted with the situation where
on both surfaces we have singular
fibers of the same type as before, but 
with pairing interchanged:
\begin{eqnarray}
\label{eq:latter}
\left.\begin{matrix}
I_1\times I_1\\
\\
I_2\times I_0
\end{matrix}
\;\;
\right\}
\;\;
III\times I_1.
\end{eqnarray}
This time the
branch of the $I_1$ singularity does not reduce to
a component of the type $III$ fiber mod $\p$.
Hence the small
resolution of the  reduction does not agree with the
reduction of the small resolution.
Consequently there is bad reduction at~$\p$.

A similar reduction pattern occurs for the merging of $I_3$ and $I_1$ to $IV$.
If paired with an unchanged semi-stable fiber, the reduction is good if and only if the original nodes in characteristic zero come from the $I_3$ fiber (see for instance \cite[p.~222]{Schuett}).

\subsection{Application to modularity}
The situation from \eqref{eq:latter} may be observed in
\cite[Sect.~3.2]{Schuett2}.
There the reduction
of some fiber product $W_2$ at the prime 359 has a
(non--projective) small resolution, but this
resolution does not agree with the reduction of the
small resolution in characteristic
zero. 
Consequently there is bad reduction at 359.
Curiously this prime does not divide the level
of the corresponding weight four cusp form, but only of the weight 2 form.
Note that the precise corresponding cusp form of weight four and level 55, as conjectured in
\cite[Conj,~3]{Schuett2}, can be verified by replacing one of the factors in the fiber product by a $2$-isogenous elliptic surface.
This construction for correspondences between Calabi-Yau threefolds was introduced in \cite{Kap}.
In the present situation, it exchanges the singular fibers of type $I_1$ and $I_2$ above;
the new fiber product is rigid of type (3) in the notation of
Proposition~\ref{Prop:Schoen}. 
Reduction modulo $359$ switches from the second degeneration type to the first,
so $359$ is a prime of good reduction for the new fiber product.
Hence the corresponding Galois representations are unramified at $359$.
By \cite{Schuett2}, the computed traces thus suffice to determine the associated cusp form.

\section{Fiber products} 
\label{s:setup}

In this section we employ the notation from \cite{Schoen}.
Let $r:Y\lra\PP$ and $r':Y'\lra\PP$ be two
rational elliptic surfaces with sections. 
Denote
by $S$ and $S'$ the images of singular fibers of
$Y$ and $Y'$ in $\PP$, and  let $S''=S\cap S'$.  
Assume that the
singular fibers of $Y$ and $Y'$ are reduced
(Kodaira types $I_n, II, III, IV$).
Then the argument from \cite{Schoen} for the semi-stable case generalize directly to this set-up to show that the fiber product
\[
X = Y \times_{\PP} Y'
\]
is simply connected and has trivial dualizing sheaf. 
Hence any crepant resolution of $X$ is a Calabi-Yau threefold.
Often one tries to work with small resolutions.
For this we fix the following

\subsection{Assumption}
\label{ass}
For
any $b\in\mathbb P^1$, either both fibers $Y_b$ and
$Y'_b$ are semistable or they have the same type.

\medskip

Under this assumption, we shall show that the fiber product
$X$ has a smooth model $\tilde
X$ which is a Calabi--Yau threefold (as a manifold
or as an algebraic space). 
Recall that the singularities of $X$ correspond to points $(x,x')$,
where $x$ (resp. $x'$) is a singular point of the
fiber $Y_b$, resp. $Y'_b$ at $b=r(x)=r'(x')$. We
shall write down a crepant resolution separately for
every type of singular fiber.
The following arguments hold true in any characteristic different from $2,3$.
In the exceptional characteristics, one might have to pay special attention, for instance 
because of the presence of wild ramification on the elliptic surfaces (cf.~\cite{SS}),
but we will not need this here.

\subsection{Fiber type $I_n\times I_m$}
Then the fiber product has $nm$ nodes on the fiber
$Y_b\times Y'_b$.
They admit a small crepant
resolution. 
Since the Euler characteristic of a
fiber of type $I_n$ is $n$, we get $e(\tilde
X_b)=2nm$. The number of irreducible components of
$\tilde X_b$ equals $nm$.

\begin{subsection}{Fiber type $II\times II$}
  The fiber $X_{b}$ has one singular point.
  The small resolution from \eqref{eq:D4} replaces the $D_4$
singularity with two intersecting lines
  so $e(\tilde X_b)=2\times2+1+1=6$.
  The 
  fiber $\tilde X_b$ is
  irreducible. (Analytically, this small
  resolution behaves like resolving $I_1\times
  IV$.)
\end{subsection}

\begin{subsection}{Fiber type $III\times III$}
 The fiber $X_{b}$ has one singular point.
 The local
(analytic) equation of $X$
near the singularity is 
\[
x^2+y^4-z^2-t^4=0.
\]
 We blow--up first the
surface
$x-z=y-t=0$.
The singular point is replaced by a line containing  two
nodes on it.
Small resolutions replace the nodes with two further lines.
Hence $e(\tilde
X_b)=3\times3+1+1+1=12$.
The fiber $\tilde X_b$
has four irreducible components.
\end{subsection}

\begin{subsection}{Fiber type $IV\times IV$}
 The fiber $X_{b}$ has one singular point, which is an ordinary
triple point.
This time we perform a big resolution.
Blowing--up
a threefold triple point yields a crepant resolution with the singular
point  replaced by a cubic in $\mathbb P^{3}$.
Hence $e(\tilde
X_b)=4\times4+8=24$.
The fiber $\tilde X_b$
has 10 irreducible components.
\end{subsection}

\subsection{Invariants of $\tilde X$}

We have seen that under the above assumption, the fiber product $X$ admits a crepant resolution $\tilde X$.
We shall now compute some invariants of the complex Calabi-Yau threefold $\tilde X$.
Our arguments follow the line of \cite{Schoen}.
We are particularly interested in the dimension of the space of infinitesimal deformations.
By Serre duality, this dimension equals $h^{1,2}(\tilde X)$.
Thus $\tilde X$ is rigid if and only if $h^{1,2}(\tilde X)=0$.

The explicit crepant resolutions of the singularities of $X$ allow us to compute the Euler characteristic
of the complex Calabi-Yau threefold $\tilde X$ explicitly as the sum of Euler characteristics of singular
fibers. 
Denote by $n_{2}$ (resp. $n_{3}$, $n_{4}$) the number of 
singular fibers of type $II$ (resp. $III, IV$). 
For $b\in\mathbb P^{1}$, let $t(b)$ (resp.
$t'(b)$) denote the number of components of the fiber $Y_{b}$
(resp. $Y'_{b}$). Then
\[
e(\tilde X)=2\left(\sum_{b\in
\PP}t(b)t'(b)-n_{2}-3n_{3}-4n_{4}\right)
=
2\left(\sum_{b\in
S''}t(b)t'(b)+2n_{2}+2n_{3}+3n_{4}\right).
\]
We shall compare this with another computation of the Euler characteristic of $\tilde X$ 
as the alternating sum of Betti numbers.
By the exponential sequence, Calabi-Yau threefolds have $b_2=h^{1,1}=\rho$
where $\rho$ denotes the Picard number.
In the present situation 
an analogue of the Shioda-Tate formula for elliptic surfaces with section 
\cite[Cor.~5.3]{ShMW}
shows 
\[
\rho(\tilde X)=d+3+\rk(\MW(Y))+\rk(\MW(Y'))+\sum_{B\in
S\cup S'} (\#(\text{components of }\tilde X_{b}) - 1).
\]
Here $d=1$ when $Y$ and $Y'$ are
isogenous (i.e.~their generic fibers are isogenous) and $d=0$ otherwise.
Under the fixed assumptions on the singular fibers,
the Shioda-Tate formula learns us that the rational elliptic surfaces $Y, Y'$ have
\[
\rk(\MW(Y))=\#S-4+n_{2}+n_{3}+n_{4},
\;\;\;
\rk(\MW(Y'))=\#S'-4+n_{2}+n_{3}+n_{4}.
\]
This gives $\rho(\tilde X)$.
On the other hand,
$b_3(\tilde X) = 2(h^{0,3}(\tilde X) + h^{1,2}(\tilde X))$ by complex conjugation 
and $h^{0,3}(\tilde X)=1$ by Serre duality.
Thus 
we deduce
\[
h^{1,2}(\tilde X) = \rho(\tilde X) - \tfrac 12 e(\tilde X).
\]
Since the number of components of the fiber $\tilde X_b$ equals $t(b)t'(b)$ unless both fibers have type IV where we have to add one for the exceptional divisor,
the latest relation can be expanded as
\begin{equation}\label{deform}
h^{1,2}(\tilde X)=d+\#(S\cup S')-5+\sum_{b\in S\setminus
S''}(t(b)-1)+\sum_{b\in S'\setminus
S''}(t'(b)-1).
\end{equation}
The above formula for the Hodge number $h^{1,2}(\tilde X)$ (under Assumption \ref{ass}) is
exactly the same as in the semistable case. 
A case-by-case analysis as in \cite{Schoen} reveals when $\tilde X$ admits no infinitesimal deformations.

\begin{prop}
\label{Prop:Schoen}
The complex Calabi--Yau threefold $\tilde X$ is rigid if and only if it
contains no fibers of type $I_{0}\times I_{n}$ or $I_{n}\times
I_{0}$ with $n>1$ and one of the following cases holds.
\begin{enumerate}
   \item $S=S'$, $\#S=4$ (then $Y$ and $Y'$ are isogenous),
   \item $\#S=\#S'=4$, $\#S''=3$ (then $Y$, $Y'$ are not
isogenous),
   \item $\#S=5$, $\#S'=4$, $S'\subset S$ (then $Y$, $Y'$ are not
isogenous),
   \item $\#S=5, S=S'=S''$, but $Y$, $Y'$ are not
isogenous,
\end{enumerate}
\end{prop}

In the sequel, we shall refer to (rigid) fiber products of type (2) and (3) according to the cases of the proposition.
Note that the isogeny conclusion of the first case fails to be valid in positive characteristic.
In fact, this failure causes non-liftability as we shall explore in the next sections.

\section{Rational elliptic surfaces with one fiber of type $II$
and three semistable fibers}
\label{s:4}

In the following sections
 we shall construct non--liftable
Calabi--Yau spaces in positive characteristic.
We achieve this by
considering fiber products of rational elliptic
surfaces which have combinations of singular
fibers which do not occur in characteristic
zero. 
Our main tool to construct these Calabi-Yau spaces are
rational elliptic surfaces with section and four singular fibers. 
Over $\mathbb C$ these surfaces have been classified by Beauville \cite{B} 
(the semistable ones, equivalently with finite Mordell-Weil group) 
and Herfurtner \cite{Herfurtner}.
In positive characteristic,
rational elliptic surfaces with finite Mordell-Weil group 
have been classified by Lang \cite{Lang-I}, \cite{Lang-II}.

In \cite{CvS}, non-liftable fiber products arising from semi-stable rational elliptic surfaces were studied.
Here we shall consider surfaces with one fiber of type $II$ and three semi-stable singular fibers.
There are four such surfaces.
M\"obius transformations (or twists) send the  fiber of type $II$ to $\infty$ and two other singular fibers to $0, 1$.
Then the remaining cusp runs through the set
\[
\left\{
\lambda, 1-\lambda, \frac 1\lambda, 1-\frac 1\lambda, \frac 1{1-\lambda}, \frac{\lambda}{\lambda-1}\right\}
\]
and the precise value determines the twist for all present surfaces.
The entries in the following tables refer to the equations in \cite{Herfurtner}.
We do not reproduce them here, since nicer models with an explicit generator of the Mordell-Weil group can be obtained from the families in Section~\ref{s:5} by specialisation, as we shall use later.

\def\arraystretch{1.4}
\begin{table}
\def\arraystretch{1.4}
\[
\begin{array}{|c|c|c|c|c|}
\hline
\parbox{1.45cm}{\centering$I_1$}&\parbox{1.45cm}{\centering$I_1$}&\parbox{1.45cm
} {\centering$I_8$}&\parbox{1.45cm}{\centering$II$}&\parbox{2.8cm}{\centering
twist}
\\
\hline1&\frac{17+56\sqrt{-2}}{81}&0&\infty&\frac{7+4\sqrt{-2}}{3t}\\
\hline0&\frac{64-56\sqrt{-2}}{81}&1&\infty&\frac{3t-(7+4\sqrt{-2})}{3t} \\
\hline0&1&\frac{8+7\sqrt{-2}}{16}&\infty&\frac{(8+7\sqrt{-2})t-27\sqrt{-2}}{
16t } \\
\hline1&0&\frac{8-7\sqrt{-2}}{16}&\infty&\frac{(8-7\sqrt{-2})t+27\sqrt{-2}}{
16t } \\
\hline\frac{64+56\sqrt{-2}}{81}&0&1&\infty&\frac{3t-(7-4\sqrt{-2})}{3t} \\
\hline\frac{17-56\sqrt{-2}}{81}&1&0&\infty&\frac{7-4\sqrt{-2}}{3t}\\\hline
\multicolumn{5}{c}{}\\
\hline
\parbox{1.45cm}{\centering$I_1$}&\parbox{1.45cm}{\centering$I_2$}&\parbox{1.45cm
} { \centering$I_7$}&\parbox{1.45cm}{\centering$II$}&\parbox{2.8cm}{\centering
twist} \\
\hline
\frac{32}{81}&1&0&\infty&-\frac{8}{9t}
\\
\hline
\frac{49}{81}&0&1&\infty&\frac{9t+8}{9t}\\\hline
1&0&\frac{81}{49}&\infty&\frac{81t+72}{49t}\\\hline
0&1&-\frac{32}{49}&\infty&\frac{-32t-72}{49t}\\\hline
0&-\frac{49}{32}&1&\infty&\frac{4t+9}{4t}\\\hline
1&\frac{81}{32}&0&\infty&-\frac{9}{4t}\\\hline
\multicolumn{5}{c}{}\\
\hline
\parbox{1.45cm}{\centering$I_1$}&\parbox{1.45cm}{\centering$I_4$}&\parbox{1.45cm
} { \centering$I_5$}&\parbox{1.45cm}{\centering$II$}&\parbox{2.8cm}{\centering
twist} \\\hline
\frac{1}{81}&1&0&\infty&\frac{1}{8t-1}\\\hline
\frac{80}{81}&0&1&\infty&\frac{8t}{8t-1}\\\hline
1&0&\frac{81}{80}&\infty&\frac{81t}{80t-10}\\\hline
0&1&-\frac{1}{80}&\infty&\frac{10+t}{10-80t}\\\hline
0&-80&1&\infty&\frac{8t+80}{8t-1}\\\hline
1&81&0&\infty&-\frac{81}{1-8t}\\\hline
\multicolumn{5}{c}{}\\
\hline
\parbox{1.45cm}{\centering$I_2$}&\parbox{1.45cm}{\centering$I_3$}&\parbox{1.45cm
} { \centering$I_5$}&\parbox{1.45cm}{\centering$II$}&\parbox{2.8cm}{\centering
twist} \\\hline
\frac{27}{32}&1&0&\infty&\frac{3}{3-t}\\\hline
\frac{5}{32}&0&1&\infty&\frac{t}{t-3}\\\hline
1&0&\frac{32}{5}&\infty&\frac{32t}{5t-15}\\\hline
0&1&-\frac{27}{5}&\infty&\frac{15+27t}{15-5t}\\\hline
0&-\frac5{27}&1&\infty&\frac{9t+5}{9t-27}\\\hline
1&\frac{32}{27}&0&\infty&\frac{32}{27-9t}\\\hline
\end{array}
\]
\caption{Rational elliptic surfaces with a $II$ fiber and three semi-stable singular fibers}
\label{tab0}
\end{table}

Let $Y, Y'$ be two different elliptic surfaces from the
above lists, determined by their loci of singular fibers
\[
S=\{0,1,\infty,\lambda\} \;\;\; \text { and } \;\;\;
S'=\{0,1,\infty,\lambda'\}.
\]
Since
$\lambda\not=\lambda'$, the elliptic  surfaces  $Y$ and $Y'$ are
not isogenous. 
By Section~\ref{s:setup}, 
the fiber product 
$X:=Y\times_{\mathbb P^1}Y'$ has a small resolution $\tilde X$ which is a
(non--projective) Calabi--Yau manifold. From 
\eqref{deform} we get
\[h^{1,2}(\tilde X)=t(\lambda)+t'(\lambda')-2,\]
so $\tilde X$ is rigid if and only if both fibers $S_\lambda$
and $S'_{\lambda'}$ have Kodaira type
$I_1$. When $\tilde X$ is defined over $\mathbb
Q$, it is modular even if
it is non--rigid, since the compatible system of Galois of $H^3(\tilde X)$ (for $\ell$-adic cohomology in the category of algebraic spaces, or for a partial big projective resolution)
can be split into two-dimensional subrepresentations by \cite{HulekVerrill}. 

The primes of bad reduction of $\tilde X$ (or a
big resolution of $X$) are primes at which either one of the factors degenerates or the
fiber product changes the configuration of
singular fibers.
The latter corresponds to the merging of some elements of the set
$\{0,1,\infty,\lambda,\lambda'\}$ upon reduction. 
These primes
can be identified as the divisors of numerator and
denominator of
$(\lambda-\lambda')\lambda\lambda'(\lambda-1)(\lambda'-1)$. 

We are interested in the situation when the
reduction $X_p$ of $\tilde X$ modulo $p$ is a
nodal Calabi--Yau space.
In our setting, this happens exactly 
when
$\lambda=\lambda'\not\in\{0,1,\infty\}$ modulo the
given prime.  
Equivalently the
prime number $p$ is a divisor of the numerator of
$\lambda-\lambda'$, but it is not a divisor of the
numerator or  denominator of
$\lambda\lambda'(\lambda-1)(\lambda'-1)$.  
In the
opposite situation when the prime $p$ divides the
numerator or the denominator of
$\lambda\lambda'(\lambda-1)(\lambda'-1)$,
we get 
certain kinds of degenerate reduction.

Table~\ref{tab1} collects the bad primes for
various fiber products of
elliptic surfaces from Table~\ref{tab0}. In the
table ``nodal bad prime'' refers to a prime divisor of the
numerator of $\lambda-\lambda'$ whereas
``degenerate bad prime'' means a prime dividing the
numerator or denominator of $\lambda\lambda'(\lambda-1)(\lambda'-1)$.
The last column gives the Hodge number
$h^{1,2}(\tilde X)$ from \eqref{deform}. 
Note that the same primes
can be obtained for several different pairs
$\lambda,\lambda'$.
In order keep the length
of the table reasonable,
 we list only
the ``simplest'' example for any given prime.

Introduce earlier:
$\alpha=\frac{17+56\sqrt{-2}}{81}$

\begin{table}
\[
{
\begin{array}{|c|c|c|c|c|}
\hline
\multirow{2}{*}{$\lambda$}&\multirow{2}{*}{$\lambda'$}&\multicolumn{2}{c|}{\text
  {Bad primes}}&\multirow{2}{*}{$h^{1,2}$}\\\cline{3-4} 
&&\text{nodal}&\text{degen.}&
\\\hline
\frac{5}{32} & \frac{32}{5} &    37 &   5,3,2&5\\\hline
\frac{5}{32} & -\frac{27}{5}&    127,7 &   5,3,2&5\\\hline
\frac{5}{32} & -\frac{5}{27} &    59 &   5,3,2&3\\\hline
\frac{5}{32} & \frac{80}{81} &    431 &   5,3,2&1\\\hline
\frac{5}{32} & 81 &    199,13 &   5,3,2&4\\\hline
\frac{5}{32} & \frac{81}{80} &    137 &   5,3,2&5\\\hline
\frac{5}{32} & \frac{1}{81} &    373 &   5,3,2&1\\\hline
\frac{5}{32} & \frac{49}{81} &    1163 &   7,5,3,2&1\\\hline
\frac{5}{32} & \frac{32}{81} &    619 &   7,5,3,2&1\\\hline
\frac{5}{32} & -\frac{32}{49} &    47 &   7,5,3,2&7\\\hline
\frac{5}{32} & \frac{81}{49} &    2347 &   7,5,3,2&7\\\hline
\frac{1}{81} &81 &     41 &   5,3,2&3\\\hline
\frac1{81} & -\frac1{80} &    23,7 &   5,3,2&4\\\hline
\frac1{81} & -{80} &    6481 &   5,3,2&3\\\hline
\frac1{81} & \frac{80}{81} &    79 &   5,3,2&0\\\hline
\frac1{81} & \frac{81}{32} &    6529 &   7,5,3,2&1\\\hline
\end{array}
\;
\begin{array}{|c|c|c|c|c|}
\hline
\multirow{2}{*}{$\lambda$}&\multirow{2}{*}{$\lambda'$}&\multicolumn{2}{c|}{\text
  {Bad primes}}&\multirow{2}{*}{$h^{1,2}$}\\\cline{3-4} 
&&\text{nodal}&\text{degen.}&
\\\hline
\frac1{81} & \frac{32}{81} &    31 &   7,5,3,2&0\\\hline
\frac1{81} & -\frac{32}{49} &    139,19 &   7,5,3,2&6\\\hline
\frac1{81} & -\frac{49}{32} &    4001 &   7,5,3,2&1\\\hline
\frac{32}{81} &\frac{81}{32} &     113 &   7,3,2&1\\\hline
\frac{32}{81} & \frac{81}{49} &    4993 &   7,3,2&6\\\hline
\frac{32}{81} & \frac{49}{81} &    17 &   7,3,2&0\\\hline
\frac{32}{81} & -\frac{32}{49} &    13, 5 &   7,3,2&6\\\hline
\alpha & \frac{32}{5} &    1187,67 &   5,3,2&4\\\hline
\alpha & \frac{32}{27} &    97,43 &   5,3,2&2\\\hline
\alpha & \frac{1}{81} &    17 &   5,3,2&0\\\hline
\alpha & -80 &    521201 &   5,3,2&3\\\hline
\alpha & \frac{80}{81} &    19,11,7 &   5,3,2&0\\\hline
\alpha & \frac{32}{81} &    89,73 &   7,5,3,2&0\\\hline
\alpha & -\frac{49}{32} &    281,107,11 &
7,5,3,2&1\\\hline
\alpha &
1-\alpha
&337&3,2&0\\\hline 
\multicolumn{4}{c}{\phantom{a}}
\end{array}}
\]
\caption{Fibre products specified by $\lambda, \lambda'$ with bad primes and $h^{1,2}$}
\label{tab1}
\end{table}

\begin{thm}
\label{Thm:tab}
Let $X$ be a fiber product as specified in Table~\ref{tab1}.
Let $p$ denote a nodal prime.
Then any small resolution of $X$ over $\bar\F_p$ is a non-liftable Calabi-Yau space.
\end{thm}

We conclude this section with a partial proof of Theorem \ref{Thm:tab}.
Let $\tilde X$ be a rigid fiber product from Table~\ref{tab1} ($h^{1,2}(\tilde X)=0$).
Let $p$ be a nodal prime.
By Theorem~\ref{thm:cvs}, the small
resolution of the reduction $\tilde X_p$ is a Calabi--Yau space in characteristic $p$ that
does not lift to characteristic 0.
Consequently we get non--liftable Calabi--Yau spaces in characteristics
$7,11,17,19,31,73,79,89,337$. 

For the non-rigid fiber products from Table~\ref{tab1},
we shall develop an alternative approach in the next section that will allow us to apply Theorem~\ref{thm:cvs} and deduce non-liftability at the nodal primes.

\section{Rigid fiber products}
\label{s:5}

For some of the fiber products from Table~\ref{tab1},
we have been able to apply Theorem~\ref{thm:cvs} 
to deduce non-liftability at the nodal primes directly.
Namely those are the cases 
where the corresponding complex Calabi-Yau threefold is already rigid ($h^{1,2}=0$).
In this section, we prove non-liftability at the nodal primes for all other fiber products from Table~\ref{tab1}.
The general idea is to lift a partial small resolution to a rigid smooth Calabi-Yau space over some number field.
Equivalently, we exhibit a
rigid Calabi-Yau threefold over some number field that reduces modulo a suitable prime to the given fiber product with an additional node.
Here we will work with fiber products of type (3) in terms of Proposition~\ref{Prop:Schoen}.
To achieve this, we shall keep one of the elliptic surfaces
of the original fiber product from Table~\ref{tab1} unchanged, whereas we replace the second one with
an elliptic surface with five singular fibers,
four semistable and one of type $II$.
We shall first introduce the rational elliptic surfaces that we will use for this construction.
Then we will explain how to complete the proof of Theorem~\ref{Thm:tab}.

\subsection{Rational elliptic surfaces with one fiber of type $II$
and four semistable singular fibers}

We shall study rational elliptic surfaces with five singular fibers under the assumption that one fiber has type $II$ and all others are semistable.
These surfaces come in seven families, each depending on one parameter.
Here we shall only introduce  four families which suffice for our purposes.

We employ the short-hand $[n_1,\hdots,n_4,II]$ for the configuration of singular fibers.
Here $n_1+\hdots+n_4=10$, since the Euler number equals $12$.
Each family will be given in terms of an extended Weierstrass form over $\Q$ with parameter $m$.
We also list the parameter choices where the family degenerates to one of the surfaces from the previous section or other interesting surfaces (see Section~\ref{s:III}).
By way of reduction, the equations stay valid outside characteristics $2$ and $3$.
Those characteristics play a special role since fibers of type $II$ do necessarily come with wild ramification by \cite{SS}.
Hence mod $2$ and $3$ there has to occur some degeneration of singular fibers.

As exploited in \cite{SS}, expanded Weierstrass forms are very useful when one wants to find elliptic surfaces with specific singular fibers.
In the present situation, the equations can be found without much difficulty,
so we decided to omit the details.

The given families have Mordell-Weil rank two.
This follows from the Shioda-Tate formula \cite[Cor.~5.3]{ShMW} since the Picard number is ten while the contribution from singular fibers and zero section amounts to 8.
From the specific Weierstrass model, it is easy to read off a section of the elliptic fibration by setting $x=0$.
This section has infinite order since fibers of type $II$ do not accommodate torsion sections of order relatively prime to the characteristic.
Unless otherwise noted, all listed degenerate members of the family share the property that the Mordell-Weil rank drops to one 
due to the fiber degeneration.
The given section specialises to a generator of the Mordell-Weil group up to 
finite index.

\subsection{[1,1,1,7,II]}
\label{ss:1,1,1,7}
Extended Weierstrass form:
\[
y^2 = x^3+(t^2+mt+3)x^2+(2t+3)x+1.
\]
\[
\begin{array}{|c|c|c|c|c|}
\hline
\parbox{2.55cm}{\centering$I_1$}&\parbox{2.55cm}{\centering$I_1$}&\parbox{2.55cm
}{\centering$I_1$}&\parbox{1.45cm
} { \centering$I_7$}&\parbox{1.45cm}{\centering$II$} \\\hline 
\multicolumn{3}{|c|}{
\parbox{7.8cm}
{\rule{0cm}{4mm}
$4mt^3-12t^3+8m^2t^2-24mt^2+3t^2+4m^3t-12m^2t+30mt-76t+27m^2-108m+108=0$
\rule[-2mm]{0cm}{1mm}
}}
&\infty&0\\
\hline
\end{array}
\]
\begin{center}
\begin{tabular}{ccccc}
\hline
$m$ && 3 & -6 & 2\\
\hline
degeneration && [1,1,8,II] & [1,2,7,II] & [1,1,7,III]\\
\hline
\end{tabular}
\end{center}

\subsection{[1,1,2,6,II]}
\label{ss:1,1,2,6}

Extended Weierstrass form:
\[
y^2 = x^3+(t^2+3t+3m)x^2+((3m-1)t+3m^2)x+m^3.
\]
\[
\begin{array}{|c|c|c|c|c|}
\hline
\parbox{3.1cm}{\centering$I_1$}&\parbox{3.1cm}{\centering$I_1$}&\parbox{1.45cm
}{\centering$I_2$}&\parbox{1.45cm
} {
  \centering$I_6$}&\parbox{1.45cm}{\centering$II$}
\\\hline
\multicolumn{2}{|c|}{
4mt^2-t^2+18mt-4t+27m^2=0}
&\frac1{m-1}\rule[-2mm]{0mm}{1mm}&\infty&0
\\\hline
\end{array}
\]
\begin{center}
\begin{tabular}{ccccc}
\hline
$m$ && 1 &1/4 & 0\\
\hline
degeneration && [1,1,8,II] & [1,2,7,II] & [1,2,6,III]\\
\hline
\end{tabular}
\end{center}

\subsection{[1,1,3,5,II]}
\label{ss:1,1,3,5}

Extended Weierstrass form
\[
y^2+(mt-27m+18)xy+16ty = x^3+tx^2.
\]
\[
\begin{array}{|c|c|c|c|c|}
\hline
\parbox{3.1cm}{\centering$I_1$}&\parbox{ 3.1cm}{\centering$I_1$}&\parbox{1.45cm
}{\centering$I_3$}&\parbox{1.45cm
} { \centering$I_5$}&\parbox{1.45cm}{\centering$II$}
\\\hline
\multicolumn{2}{|c|}{
\parbox{5cm}
{$m^4t^2
-2m^2(27m^2-28m-4)t
+(27m-2)(3m-2)^3=0$}}
&0&\infty&27
\\\hline
\end{array}
\]
\begin{center}
\begin{tabular}{ccccc}
\hline
$m$ && 2/27 & 1 & 1/9\\
\hline
degeneration && [1,4,5,II] & [2,3,5,II] & [1,3,5,III]\\
\hline
\end{tabular}
\end{center}

\subsection{[1,2,3,4,II]}
\label{ss:1,2,3,4}

Extended Weierstrass form
\[
y^2+(mt(m+8)+1-m)xy+\tfrac{27}4m^2t(t(m+8)-1)y =
x^3+\tfrac{27}4tmx^2.
\]
\[
\begin{array}{|c|c|c|c|c|c|}
\hline
\parbox{2.45cm}{\centering$I_1$}&\parbox{ 1cm}{\centering$I_2$}&\parbox{1.cm
}{\centering$I_3$}&\parbox{1.cm
} { \centering$I_4$}&\parbox{1.cm}{\centering$II$}&\parbox{1.2cm}{\centering
twist} \\\hline
\rule[-3mm]{0cm}{8mm}\parbox{2.3cm}{$\frac{4(m-1)^3}{m(4m+5)(m+8)^2} $}
&\frac1{m+8}&0&\infty &\frac1m&\\\hline
\rule[-3mm]{0cm}{8mm}-\frac{4(m-1)^3}{81(m+2)^2}&-\frac m8&0&1&\infty&\frac{mt}{mt-1}
\\\hline
\end{array}
\]
\begin{center}
\begin{tabular}{ccccccc}
\hline
$m$ && 0 & -5/4 & -2 & 1 & -2/7\\
\hline
degeneration && [1,4,5,II] & [2,3,5,II] & [2,3,4,III] & [2,4,II,IV] & [3,4,II,III] 
\\
\hline
\end{tabular}
\end{center}

Note that the last two specialisations do not cause a drop of the Mordell-Weil rank.
For the specialisation at $m=0$, we work with the twisted Weierstrass form.
Here one can eliminate factors of $m$ to obtain the following model with nodal specialisation at $m=0$:
\[
y^2+4(9t+m-1)xy+27t(t-1)(8t+m)y = x^3 -3(72t^2+9(2m-1)t+(m-1)^2) x^2.
\]

\subsection{Illustration of the method}
\label{ss:illustration}

We shall construct fiber products of type (3) in the notation of Proposition~\ref{Prop:Schoen} from
\begin{itemize}
\item 
$Y$: the family of rational elliptic surfaces of type
[1,2,3,4,II],
\item
$Y'$: the elliptic surfaces with four singular fibers
from the previous section.
\end{itemize}
Working with the twisted model of the former family,
we have 
\[
S=\left\{0,1,\infty,-\frac
m8,-\frac{4(m-1)^{3}}{81(m+2)^{2}}\right\}
\;\;\; \text{ and } \;\;\;
S':=\{0,1,\infty,\lambda\}.
\]
After a further M\"obius transformation for the second elliptic surface, 
we can assume that the fiber of type $II$ sits at $\infty$
for both elliptic surfaces $Y, Y'$ (as in the previous section).

We now
choose $m=-8\lambda$.
Then
the fiber product $X:=Y\times_{\mathbb P^{1}}Y'$ has a small
resolution $\tilde X$ which is a non--projective rigid
Calabi--Yau manifold by Proposition~\ref{Prop:Schoen}. 
In order to find the primes of bad reductions we have to find the primes at
which the construction degenerates.
In the ring of integers of the field of
definition of the algebraic space $\tilde X$, the candidates for bad reduction are exactly the
prime divisors of the numerators and denominators of
the numbers 
\[
-\frac m8, -\frac m8-1,
-\frac{4(m-1)^{3}}{81(m+2)^{2}},
-\frac{4(m-1)^{3}}{81(m+2)^{2}}-1,
-\frac{4(m-1)^{3}}{81(m+2)^{2}}+\frac m8.
\]
For the reduction we distinguish four cases:
\begin{enumerate}
\item
If $-\frac m8$ equals $0, 1$ or $\infty$ modulo a given prime $\p$, 
then $Y'$ reduces to an elliptic surface with only three singular fibers.
Then either there is a non-reduced singular fiber 
not allowing a small resolution (cf.~Assumption \ref{ass})
or the characteristic is at most $7$
(with configurations such as [7,II,III] in characteristic 7 
or [5,5,II], [5,III,IV] in characteristic 5). 
In particular,
we always get a bad prime, but we cannot get a non--liftable
Calabi--Yau space in characteristic greater then 7.
\item
If $-\frac{4(m-1)^3}{81(m+2)^2}=1$, i.e.~$m=\frac 54$,
the surface $Y$ degenerates to
[2,3,5,II].
This imposes an additional node on the fiber product.
By Theorem~\ref{thm:cvs}, reduction gives a Calabi--Yau space in positive
characteristic $p$ that does not lift to characteristic zero 
(one of the examples
predicted in Table~\ref{tab1}).

\item
If $-\frac{4(m-1)^3}{81(m+2)^2}=\infty$, i.e.~$m=-2$, 
then the elliptic surface $Y$ degenerates to
[2,3,4,III]. 
Since a fiber of type $II\times III$ has no small resolution by \cite{Kap0}, 
$\p$ is a bad prime, but we cannot apply Theorem~\ref{thm:cvs}.

\item
If $-\frac{4(m-1)^3}{81(m+2)^2}=-\frac m8$ (resp.
$-\frac{4(m-1)^3}{81(m+2)^2}=0$), 
the surface $Y$ degenerates to
[III,3,4,II] (resp. [2,IV,4,II]). 
By Section~\ref{ss:359}, the fiber of type $I_n\times III$ (resp.
$I_n\times IV$) admits a small resolution 
in a way that is compatible with the small resolution of the original fiber of type 
$I_n\times I_2$ (resp. $I_n\times I_3$).
In this case, $\p$ is a prime of good reduction.
\end{enumerate}
 
Table~\ref{tab2} list choices of $\lambda, m$ such that the fiber product of $Y$ and $Y'$ is rigid (type (3)).
We then give nodal and degenerate primes.
 By our previous considerations, we obtain non--projective rigid Calabi--Yau spaces with all bad primes from both columns and non--liftable Calabi--Yau spaces in all nodal characteristics.
 
\begin{table}
\[
{
\def\arraystretch{1.5}
\begin{array}{|c|c|l|l|l|}
 \hline\multirow{2}{*}{\parbox{2cm}{\centering
$\lambda$}}&\multirow{2}{*}{\parbox{2cm}{\centering
$m$}}&\multicolumn{2}{|c|}{\text{bad primes}}\\\cline{3-4}
&&\parbox{2cm}{nodal}&\parbox{2cm}{degenerate}\\
\hline
-\frac{27}{5} & \phantom{-}\frac{216}{5}  & 127,7& 113,5,3,2 \\\hline
\phantom{-}\frac{32}{5} & -\frac{256}{5} & 37& 41,5,3,2  \\\hline
-\frac{5}{27} & \phantom{-}\frac{40}{27} & 59& 47,5,3,2   \\\hline
\phantom{-}\frac{32}{27} & -\frac{256}{27}  & 127,7& 101,5,3,2  \\\hline
\phantom{-}\frac{1}{81} & -\frac{8}{81} & 373& 11,7,5,3,2  \\\hline
\phantom{-}81 & -648 & 199,13  & 19,17,5,3,2 \\\hline
\phantom{-}\frac{81}{80} & -\frac{81}{10} & 137& 61,5,3,2   \\\hline
\phantom{-}\frac{80}{81} & -\frac{640}{81} & 431& 239,5,3,2  \\\hline
\phantom{-}\frac{32}{81} & -\frac{256}{81} & 619& 47,7,3,2   \\\hline
\phantom{-}\frac{49}{81} & -\frac{392}{81} & 1163& 23,7,5,3,2   \\\hline
\phantom{-}\frac{81}{32} & -\frac{81}{4} & 19& 73,7,3,2   \\\hline
\phantom{-}\frac{81}{49} & -\frac{648}{49} & 2347& 11,7,5,3,2   \\\hline
-\frac{32}{49} & \phantom{-}\frac{256}{49} & 47& 59,7,3,2   \\\hline
\end{array}}
\]
\caption{Rigid fiber products involving [1,2,3,4,II]}
\label{tab2}
\end{table}

\begin{rem}
Although the singular Calabi--Yau threefolds from Table~\ref{tab2} are defined over $\mathbb
Q$, as well as the big resolution (non-crepant), it is a delicate question
if the small resolution can be defined over $\mathbb Q$.
For modularity questions, however, it suffices to consider the (partial) big resolution, a smooth projective variety with $b_3=2$.
Modularity then follows from Serre's conjecture \cite{KW}.
We found the big bad primes of the threefolds  remarkable.
Necessarily they appear in the level of the associated modular form.
In comparison, it seems that for all examples so far the biggest bad prime (up to twists) was $73$ (cf.~\cite{Schuett}).
\end{rem}

\subsection{Proof of Theorem~\ref{Thm:tab}}

We shall now show how to set up other rigid fiber products in order to prove the non-liftability of all nodal examples from the first table and thus complete the proof of Theorem~\ref{Thm:tab}.

Let $Y, Y'$ denote the rational elliptic surfaces giving a fiber product $X$ from Table~\ref{tab1} specified by $\lambda,\lambda'$. 
Assume that $h^{1,2}(\tilde X)>0$ and that $p$ denotes a nodal prime.
Then reduction mod $p$ imposes additional nodes on $\tilde X_p$, but Theorem~\ref{thm:cvs} does not apply since the original Calabi-Yau threefold $\tilde X$ is not rigid.
To prove non-liftability, we shall extend the construction from~\ref{ss:illustration}.

\begin{prop}
\label{prop}
Let $\tilde X$ be as above.
Then there is a rigid fiber product of type (3) over some number field
which reduces to $\tilde X_p$ modulo some prime $\p$ above $p$.
\end{prop}

As a corollary of the proposition, we deduce from Theorem~\ref{thm:cvs} that the reduction of $\tilde X$ mod $\p$ has a non-liftable small resolution.
This completes the proof of Theorem~\ref{Thm:tab}.

The main idea to prove Proposition~\ref{prop} is to replace $Y$ or $Y'$ in the fiber product $X$ by a suitable member of one of the families from~\ref{ss:1,1,1,7}--\ref{ss:1,2,3,4}, as we did in~\ref{ss:illustration}.
First we need the following observation:

\begin{lem}
\label{lem:spec}
Each rational elliptic surface from Section~\ref{s:4} arises from one of the families in~\ref{ss:1,1,1,7}--\ref{ss:1,2,3,4} by degeneration where an $I_1$ fiber is merged with another fiber.
\end{lem}

\begin{proof}[Proof of Proposition~\ref{prop}]
Since we only consider fiber products with $h^{1,2}(\tilde X)>0$, 
an inspection of Table \ref{tab1} reveals that
we can assume that $Y$ does not have configuration [1,1,8,II]
after exchanging $Y, Y'$ if necessary.
This will be convenient because then all singular fibers of $Y$ sit above rational points of $\PP$,
i.e.~$\lambda\in\Q$.

Let $\mathcal{Y}(m)$ be a family degenerating to $Y'$ as in Lemma~\ref{lem:spec}, say $Y'=\mathcal{Y}(m_0)$.
By choosing an appropriate $m$, 
we want to match the singular fibers of $Y$ and $\mathcal Y(m)$ in such a way 
that their fiber product has type (3) 
and the non-degenerate singular fibers for the specialisation at $m_0$ agree with those of $\tilde X$.
For the [1,2,3,4,II] family, this was easily achieved in~\ref{ss:illustration} since all cusps are rational.
For the other families, we have to be a little more careful.

Let $\mathcal Y(m)$ be one of the families from~\ref{ss:1,1,2,6} and~\ref{ss:1,1,3,5}.
A M\"obius transformation sends the reducible fibers to $0$ and $1$ and the fiber of type $II$ to $\infty$ (accordingly for $Y$ and $Y'$).
This locates the two fibers of type $I_1$  at the zeroes of the quadratic polynomial
$\Delta(t,m)$ with coefficients in $\Q(m)$.
For $Y\times_{\PP} \mathcal Y(m)$ to have type (3), we require one of these zeroes to be $\lambda$.
In other words, $m$ has to solve the equation
\begin{eqnarray}
\label{eq:t2}
\Delta(\lambda, m) = 0.
\end{eqnarray}
By reduction, we can also consider (\ref{eq:t2}) as an equation over $\F_p$.
Since $\lambda\equiv\lambda'\mod p$, the reduced equation has the solution $m_0$ mod $p$.
This solution lifts to a solution $m_1$ of (\ref{eq:t2}) over some number field.
Then, by construction, the fiber product $Y\times_{\PP} \mathcal Y(m_1)$ has type (3) and reduces to $X$ modulo some prime above $p$. 

Finally we consider the case where $\mathcal Y(m)$ has type [1,1,1,7,II].
As before we can take the singular fibers of types $I_7$ and $II$ to $0$ and $\infty$ by a M\"obius transformation.
Then the three $I_1$ fibers of $\mathcal Y(m)$ are located at the zeroes of the cubic polynomial 
$\Delta(t,m)$ with coefficients in $\Q(m)$.
In order to set up a fiber product with $Y$ of type (3), we need that two of these zeroes are multiples by $\lambda$.
This condition is encoded in the product of irreducible factors of the
resultant $\mbox{Res}(\Delta(t,m),\Delta(\lambda t,m))$ with respect
to $t$ which do not divide $\Delta(0,m)$. Equivalently for any
sufficiently large $N$
\[
F(\lambda,m) = \mbox{Res}(\Delta(t,m),\Delta(\lambda t,m))/(\text{gcd of resultant and }\Delta(0,m)^N).
\]
We continue as before by considering this polynomial equation over $\F_p$ and lifting the solution $m_0$ mod $p$ to some number field.
The resulting fiber product has the required properties.
\end{proof}

\begin{rem}
The above degeneration technique was not required in \cite{CvS} for the following reason: 
For a semi-stable rational elliptic surface with finite Mordell-Weil group, 
one can construct an isogenous surface with an $I_1$ fiber at any given position.
Hence one can replace fiber products of such rational elliptic surfaces of type (2) in the notation of Proposition~\ref{Prop:Schoen}
by an isogenous rigid fiber product.
\end{rem}

\section{Proof of Theorem~\ref{thm}}
\label{s:III}

Theorem~\ref{Thm:tab} implies the second statement of Theorem~\ref{thm} about Table~\ref{tab1}.
For the first statement, only the primes $2, 3, 29, 61, 71$ are missing from the list of nodal primes in Table~\ref{tab1}.
As explained before, it is known that there are non-liftable
projective Calabi-Yau threefolds in characteristics $2, 3$
(cf.~\cite{Hirokado,Schroeer}).
It remains to prove the corresponding fact for Calabi-Yau spaces in characteristics $29, 61, 71$.

To achieve this, we pursue the same approach as in Section~\ref{s:4}, but this time for rational elliptic surfaces with one singular fiber of type $III$.
We shall use the specialisations of the families in~\ref{ss:1,1,1,7} --~\ref{ss:1,1,3,5} to construct fiber products of type (2).

In detail, we find the following specialisations with cusps normalised.
Here $\omega=((1-\sqrt{-7})/(1+\sqrt{-7}))^7$. 
\begin{table}[ht!]
\[\def\arraystretch{1.4}
\begin{array}{|c|c|c|c|c|}
\hline
\parbox{1.45cm}{\centering$I_1$}&\parbox{1.45cm}{\centering$I_1$}&\parbox{1.45cm
} { \centering$I_7$}&\parbox{1.45cm}{\centering$III$}&\parbox{2.8cm}{\centering
twist} \\\hline
\omega &1&0&\infty&\\\hline
1-\omega&0&1&\infty&1-t\\\hline
\multicolumn{5}{c}{}\\
\hline
\parbox{1.45cm}{\centering$I_1$}&\parbox{1.45cm}{\centering$I_2$}&\parbox{1.45cm
} { \centering$I_6$}&\parbox{1.45cm}{\centering$III$}&\parbox{2.8cm}{\centering
twist} \\\hline
1/4&1&0&\infty&\\\hline
3/4&0&1&\infty&1-t\\\hline
\multicolumn{5}{c}{}\\
\hline
\parbox{1.45cm}{\centering$I_1$}&\parbox{1.45cm}{\centering$I_3$}&\parbox{1.45cm
} { \centering$I_5$}&\parbox{1.45cm}{\centering$III$}&\parbox{2.8cm}{\centering
twist} \\\hline
3/128 &1&0&\infty&\\\hline
125/128 &0&1&\infty&1-t\\\hline
\end{array}
\]
\caption{Rational elliptic surfaces with a $III$ fiber and three semi-stable singular fibers}
\label{tabIII}
\end{table}
We obtain several interesting fiber products of type (2). 
As before they are specified by $\lambda, \lambda'$ and yield the following nodal primes:

\begin{table}[ht!]
\[
{
\def\arraystretch{1.4}
\begin{array}{|c|c|c|c|}
\hline
\multirow{2}{*}{$\lambda$}&\multirow{2}{*}{$\lambda'$}&\multicolumn{2}{c|}{\text
  {Bad primes}}\\\cline{3-4} 
&&\text{nodal}&\text{degenerate}
\\\hline
1/4 & 3/128 &    29 &   5,3,2\\\hline
3/4 & 3/128 &    31 &   5,3,2\\\hline
3/128 & 125/128 &    61 &   5,3,2\\\hline
\omega & 1/4 & 631 & 7,5,3,2\\\hline
\omega & 3/4 & 1061 & 7,5,3,2\\\hline
\omega & 3/128 & 757, 11 & 7,3,2\\\hline
\omega & 125/128 & 151,71 & 7,3,2\\\hline
\end{array}}
\]
\caption{Rigid fiber products with $III\times III$ specified by $\lambda, \lambda'$}
\label{tab3}
\end{table}

By Theorem~\ref{thm:cvs}
the nodal primes give non-liftable fiber products.
This completes the proof of Theorem~\ref{thm}.

We conclude this paper by noting that our techniques  can be applied to other rational elliptic surfaces with four singular fibers.
For instance, one could consider the remaining twisted fiber products of the above surfaces.
This would give further nodal primes $131, 137, 503, 509, 1019,  2287, 2671$.
However, our approach does not result in any other small primes ($<100)$ or big primes ($>7000$).
Hence we decided not to include these computations in this paper.

\subsection*{Acknowledgement}

We thank Chad Schoen and Duco van Straten for many useful discussions.
We are grateful to the referee for helpful comments and suggestions.
During the preparation of this paper, the second author held positions at Harvard University 
and University of Copenhagen.
He also enjoyed the hospitality of  Jagiellonian University  in June 2008.

\end{document}